\documentclass{article}
\usepackage{amssymb}
\usepackage{amsmath}
\usepackage{graphicx}

\input{tcilatex}
\begin{document}

\begin{center}
\ \ \ \ \ \ \ \ \ \ \ \ \ \ \ \ \ \ \ \ \ \ \ \ \ \ \ \ \ \ \ \ \ \ \ \ \ \
\ \ \ \ \ \ \ \ \ \ \ \ \ 

\bigskip\ \ \ \textbf{Veli B. Shakhmurov}\ \ \ 

Okan University, Department of Mechanical Engineering, Akfirat, Tuzla 34959
Istanbul, Turkey, E-mail: veli.sahmurov@okan.edu.tr

\bigskip \textbf{Regularity properties of singular degenerate abstract
differential equations and applications\ }

\ \ \ \ \ \ \ \ \ \ \ \ \ \ \ \ \ \ \ \ 

\ \ \ \ \ \ \ \ \ \ \ \ \ \ \ \ \ \ \ \ \ \ \ 

\textbf{AMS: 34G10, 35J25, 35J70\ \ }\ \ \ \ \ \ \ \ \ \ \ \ \ 

\textbf{Abstract}
\end{center}

Singular degenerate differential operator equations are studied. The uniform
separability of boundary value problems for degenerate elliptic equation and
optimal regularity properties of Cauchy problem for degenerate parabolic
equation are obtained. These problems have a numeros applications which
occur in fluid mechanics, environmental engineering and atmospheric
dispersion of pollutants.


\begin{center}
\ \ \textbf{1. Introduction, notations and background }
\end{center}

In this work, boundary value problems (BVPs) for singular degenerate
elliptic differential-operator equations (DOEs) and the Cauchy problem for
degenerate abstract parabolic equation\ are considered. BVPs for DOEs have
been studied extensively by many researchers (see e.g. $\left[ 1-20\right] $
and the references therein). A comprehensive introduction to the DOEs and
historical references may be found in $\left[ 4\right] $ and $\left[ 6\right]
.$ The maximal regularity properties for differential operator equations
have been investigated e.g. in $\left[ 3-4,\text{ }8-17\right] $. The main
objective of the present paper is to discuss BVPs for the following singular
degenerate DOE 
\begin{equation}
-\dsum\limits_{k=1}^{n}x_{k}^{2\alpha _{k}}\frac{\partial ^{2}u}{\partial
x_{k}^{2}}+x_{k}^{\alpha _{k}}A_{k}\left( x\right) \frac{\partial u}{%
\partial x_{k}}+Au=f\left( x\right) ,  \tag{1.1}
\end{equation}%
where $A$, $A_{k}$ are linear operators in a Banach space $E.$

Several conditions for the unifororm separability and the resolvent estimate
in $E$-valued $L_{p}$-spaces are given. Especially, it is proven that
corresponding differential operator is $R$-positive and\ also is a generator
of the analytic semigroup.

By using separability properties of elliptic problem $\left( 1.1\right) $,
maximal regularity properties of Cauchy problem is derived for the singular
degenerate parabolic equation 
\begin{equation}
\frac{\partial u}{\partial t}-\dsum\limits_{k=1}^{n}x_{k}^{2\alpha _{k}}%
\frac{\partial ^{2}u}{\partial x_{k}^{2}}+Au=f\left( t,x\right) .  \tag{1.2}
\end{equation}%
One of the important characteristics of these DOEs considered here are that
the degeneration process take place at different speeds on boundaries, in
general. Unlike the regular degenerate equations, due to singularity of
degeneration, boundary conditions only on undegenerate points are given.
Therefore, only one boundary condition with respect to the given variable is
given. Note that, maximal regularity properties for nonlinear DOEs are
studied e.g. in $\left[ 1,7,9-11,15\right] .$

In application, the BVP for infinity system of singular degenerate partial
differential equations and Wentzell-Robin type BVP for singular degenerate
partial differential equations on cylindrical domain are studied.

Let we choose $E=L_{2}\left( 0,1\right) $ and $A$ to be differential
operator with generalized Wentzell-Robin boundary condition defined by 
\begin{equation*}
D\left( A\right) =\left\{ u\in W_{2}^{2}\left( 0,1\right) ,\text{ }%
B_{j}u=Au\left( j\right) +\dsum\limits_{i=0}^{1}\alpha _{ij}u^{\left(
i\right) }\left( j\right) =0,\text{ }j=0,1\right\} ,\text{ }
\end{equation*}%
\begin{equation*}
\text{ }A\left( x\right) u=a\left( x,y\right) u^{\left( 2\right) }+b\left(
x,y\right) u^{\left( 1\right) }+c\left( x,y\right) u,\text{ for all }x\in
R^{n},
\end{equation*}%
where $\alpha _{ij}$ are complex numbers, $a\left( x,.\right) ,$ $b\left(
x,.\right) $, $c\left( x,.\right) $ are complex-valued functions on $\left(
0,1\right) $ for all $x\in R^{n}$. Then, from $\left( 1.2\right) $\ we get
the following Wentzell-Robin type mixed problem for singulyar degenerate
parabolic equation

\begin{equation}
\frac{\partial u}{\partial t}-\dsum\limits_{k=1}^{n}x_{k}^{2\alpha _{k}}%
\frac{\partial ^{2}u}{\partial x_{k}^{2}}+a\frac{\partial ^{2}u}{\partial
y^{2}}+b\frac{\partial u}{\partial y}+cu=f\left( x,y,t\right) ,\text{ } 
\tag{1.3}
\end{equation}%
\ \ \ 

\begin{equation*}
L_{k}u=0,\text{ }x\in G,\text{ }y\in \left( 0,1\right) ,\text{ }t\in \left(
0,T\right)
\end{equation*}%
\begin{equation}
B_{j}u=0\text{, }j=0,1,\text{ for }x\in R^{n},\text{ }  \tag{1.4}
\end{equation}%
where $L_{k}$ are boundary conditions with respect $x\in G\subset R^{n}$
that will be definet in late. By virtue of Theorem 3.1 derived here, we
obtain that problem $\left( 1.3\right) -\left( 1.4\right) $ is maximal
regular in $L_{\mathbf{p}}\left( \tilde{\Omega}\right) $, where $L_{\mathbf{p%
}}\left( \tilde{\Omega}\right) $ denotes the space of all $\mathbf{p}$%
-summable complex-valued\ functions with mixed norm and \ 
\begin{equation*}
\tilde{\Omega}=G\times \left( 0,T\right) \times \left( 0,1\right) ,\text{ }%
\mathbf{\tilde{p}=}\left( p,p_{1},2\right) .
\end{equation*}

Note that, the regularity properties of Wentzell-Robin type BVP for elliptic
equations were studied e.g. in $\left[ \text{24, 25}\right] $ and the
references therein.

Let $\gamma =\gamma \left( x\right) $ be a positive measurable function on a
domain $\Omega \subset R^{n}.$ Here, $L_{p,\gamma }\left( \Omega ;E\right) $
denote the space of strongly measurable $E$-valued functions that are
defined on $\Omega $ with the norm

\begin{equation*}
\left\Vert f\right\Vert _{L_{p,\gamma }}=\left\Vert f\right\Vert
_{L_{p,\gamma }\left( \Omega ;E\right) }=\left( \int \left\Vert f\left(
x\right) \right\Vert _{E}^{p}\gamma \left( x\right) dx\right) ^{\frac{1}{p}},%
\text{ }1\leq p<\infty .
\end{equation*}

For $\gamma \left( x\right) \equiv 1$ the space $L_{p,\gamma }\left( \Omega
;E\right) $ will be denoted by $L_{p}=L_{p}\left( \Omega ;E\right) .$

\ The Banach space\ $E$ is called an $UMD$-space if\ the Hilbert operator $%
\left( Hf\right) \left( x\right) =\lim\limits_{\varepsilon \rightarrow
0}\int\limits_{\left\vert x-y\right\vert >\varepsilon }\frac{f\left(
y\right) }{x-y}dy$ \ is bounded in $L_{p}\left( R,E\right) ,$ $p\in \left(
1,\infty \right) $ ( see. e.g. $\left[ 21\right] $ ). $UMD$ spaces include
e.g. $L_{p}$, $l_{p}$ spaces and Lorentz spaces $L_{pq},$ $p$, $q\in \left(
1,\infty \right) $.

Let $\mathbb{C}$ be the set of the complex numbers and\ 
\begin{equation*}
S_{\varphi }=\left\{ \lambda ;\text{ \ }\lambda \in \mathbb{C}\text{, }%
\left\vert \arg \lambda \right\vert \leq \varphi \right\} \cup \left\{
0\right\} ,0\leq \varphi <\pi .
\end{equation*}

\ Let $E_{1}$ and $E_{2}$ be two Banach spaces. $L\left( E_{1},E_{2}\right) $
denotes the space of bounded linear operators from $E_{1}$ into $E_{2}.$ For 
$E_{1}=E_{2}=E$ it will be denoted by $L\left( E\right) .$

A linear operator\ $A$ is said to be $\varphi $-positive in a Banach\ space $%
E$ with bound $M>0$ if $D\left( A\right) $ is dense on $E$ and $\left\Vert
\left( A+\lambda I\right) ^{-1}\right\Vert _{L\left( E\right) }\leq M\left(
1+\left\vert \lambda \right\vert \right) ^{-1}$ $\ $for any $\lambda \in
S_{\varphi },$ $0\leq \varphi <\pi ,$ where $I$ is the identity operator in $%
E$. Sometimes $A+\lambda I$\ will be written as $A+\lambda $ and will be
denoted by $A_{\lambda }$. It is known $\left[ \text{22, \S 1.15.1}\right] $
that a positive operator $A$ has well-defined fractional powers\ $A^{\theta
}.$

Let $E\left( A^{\theta }\right) $ denote the space $D\left( A^{\theta
}\right) $ with norm 
\begin{equation*}
\left\Vert u\right\Vert _{E\left( A^{\theta }\right) }=\left( \left\Vert
u\right\Vert ^{p}+\left\Vert A^{\theta }u\right\Vert ^{p}\right) ^{\frac{1}{p%
}},1\leq p<\infty ,\text{ }0<\theta <\infty .
\end{equation*}

Let $E_{1}$ and $E_{2}$ be two Banach spaces. Now $\left( E_{1},E_{2}\right)
_{\theta ,p}$, $0<\theta <1,1\leq p\leq \infty $ will denote interpolation
spaces obtained from $\left\{ E_{1},E_{2}\right\} $ by the $K$ method \ $%
\left[ \text{22 ,\S 1.3.1}\right] $.

Let $\mathbb{N}$ denote the set of natural numbers and $\left\{
r_{j}\right\} $ is a sequence of independent symmetric $\left\{ -1,1\right\} 
$-valued random variables on $\left[ 0,1\right] $. A set $K\subset L\left(
E_{1},E_{2}\right) $ is called $R$-bounded if there is a positive constant $%
C $ such that for all $T_{1},T_{2},...,T_{m}\in K$ and $%
u_{1,}u_{2},...,u_{m}\in E_{1},$ $m\in \mathbb{N}$ 
\begin{equation*}
\int\limits_{0}^{1}\left\Vert \sum\limits_{j=1}^{m}r_{j}\left( y\right)
T_{j}u_{j}\right\Vert _{E_{2}}dy\leq C\int\limits_{0}^{1}\left\Vert
\sum\limits_{j=1}^{m}r_{j}\left( y\right) u_{j}\right\Vert _{E_{1}}dy.
\end{equation*}

The smallest $C$ for which the above estimate holds is called a $R$-bound of
the collection $K$ and denoted by $R\left( K\right) .$

The $\varphi $-positive operator $A$ is said to be $R$-positive in a Banach
space $E$ if the set $L_{A}=\left\{ \xi \left( A+\xi I\right) ^{-1}\text{: }%
\xi \in S_{\varphi }\right\} ,$ $0\leq \varphi <\pi $ is $R$-bounded.

Assume $E_{0}$ and $E$ are two Banach spaces and $E_{0}$ is continuously and
densely embedded into $E$ and $\Omega $ is a domain in $R^{n}$. Let $\alpha
_{k}=\alpha _{k}\left( x\right) $ be a positive measurable functions on $%
\Omega $. Consider the Sobolev-Lions type space\ $W_{p,\alpha }^{m}\left(
\Omega ;E_{0},E\right) ,$ consisting of all functions $u\in L_{p}\left(
\Omega ;E_{0}\right) $ that have generalized derivatives $D_{k}^{m}u=\frac{%
\partial ^{m}u}{\partial x_{k}^{m}}\in L_{p,\alpha _{k}}\left( \Omega
;E\right) $ with the norm 
\begin{equation*}
\ \left\Vert u\right\Vert _{W_{p,\alpha }^{m}\left( \Omega ;E_{0},E\right)
}=\left\Vert u\right\Vert _{L_{p}\left( \Omega ;E_{0}\right)
}+\dsum\limits_{k=1}^{n}\left\Vert \alpha _{k}^{m}\frac{\partial ^{m}u}{%
\partial x_{k}^{m}}\right\Vert _{L_{p}\left( \Omega ;E\right) }<\infty .
\end{equation*}

Let $\chi =\chi \left( t\right) $ be a positive measurable function on $%
\left( 0,a\right) $ and 
\begin{equation*}
u^{\left[ i\right] }\left( t\right) =\left( \chi \left( t\right) \frac{d}{dt}%
\right) ^{i}u\left( t\right) .
\end{equation*}
Consider the following weighted abstract space 
\begin{equation*}
W_{p,\chi }^{\left[ m\right] }\left( 0,a;E_{0},E\right) =\left\{ u;u\in
L_{p}\left( 0,a;E_{0}\right) \right. ,\ u^{\left[ m\right] }\in L_{p}\left(
0,a;E\right) ,
\end{equation*}

\begin{equation*}
\left\Vert u\right\Vert _{W_{p,\chi }^{\left[ m\right] }}=\left. \left\Vert
u\right\Vert _{L_{p}\left( 0,a;E_{0}\right) }+\left\Vert u^{\left[ m\right]
}\right\Vert _{L_{p}\left( 0,a;E\right) }<\infty \right\} .
\end{equation*}

\begin{equation*}
W_{p,\chi }^{m}\left( 0,a;E_{0},E\right) =\left\{ u;u\in L_{p,\chi }\left(
0,a;E_{0}\right) \right. ,\ u^{\left( m\right) }\in L_{p,\chi }\left(
0,a;E\right) ,
\end{equation*}

\begin{equation*}
\left\Vert u\right\Vert _{W_{p,\chi }^{m}}=\left. \left\Vert u\right\Vert
_{L_{p,\chi }\left( 0,a;E_{0}\right) }+\left\Vert u^{\left( m\right)
}\right\Vert _{L_{p,\chi }\left( 0,a;E\right) }<\infty \right\} .
\end{equation*}

Let 
\begin{equation*}
\alpha =\left( \alpha _{1},\alpha _{2},...,\alpha _{n}\right) \text{, }D^{%
\left[ \alpha \right] }=D_{1}^{\left[ \alpha _{1}\right] }D_{2}^{\left[
\alpha _{2}\right] }...D_{n}^{\left[ \alpha _{n}\right] },\text{ }D_{k}^{%
\left[ i\right] }.=\left( \gamma _{k}\left( x\right) \frac{\partial }{%
\partial x_{k}}\right) ^{i}.
\end{equation*}

Consider the space\ $W_{p,\gamma }^{\left[ m\right] }\left( \Omega
;E_{0},E\right) ,$ consisting of all functions $u\in L_{p}\left( \Omega
;E_{0}\right) $ that have generalized derivatives $D_{k}^{\left[ m\right]
}u\in L_{p,}\left( \Omega ;E\right) $ with the norm 
\begin{equation*}
\ \left\Vert u\right\Vert _{W_{p,\gamma }^{\left[ m\right] }\left( \Omega
;E_{0},E\right) }=\left\Vert u\right\Vert _{L_{p}\left( \Omega ;E_{0}\right)
}+\dsum\limits_{k=1}^{n}\left\Vert D_{k}^{\left[ m\right] }u\right\Vert
_{L_{p}\left( \Omega ;E\right) }<\infty .
\end{equation*}

By reasoning as $\left[ \text{13, Theorem 2.3}\right] $ we obtain

\textbf{Theorem B}. Suppose the following conditions are satisfied:

(1) $E$ is an UMD space and\ $A$ is an $R$-positive operator in $E;$

(3)\ $\gamma =\left( \gamma _{1},\gamma _{2},...,\gamma _{n}\right) ,$ $%
\gamma _{k}\left( x\right) =\left\vert x_{k}\right\vert ^{\nu _{k}}$, $\nu
_{k}>1$ and $m$ is an integer, $\varkappa =\frac{\left\vert \alpha
\right\vert }{m}\leq 1,$ $1<p<\infty $;

(4)\ $\Omega \subset R^{n}$ is a region such that there exists a bounded
linear extension operator from $W_{p,\gamma }^{\left[ m\right] }\left(
\Omega ;E\left( A\right) ,E\right) $ to $W_{p,\gamma }^{\left[ m\right]
}\left( R^{n};E\left( A\right) ,E\right) $.

Then, the embedding $D^{\left[ \alpha \right] }W_{p,\gamma }^{\left[ m\right]
}\left( \Omega ;E\left( A\right) ,E\right) \subset L_{p}\left( \Omega
;E\left( A^{1-\varkappa -\mu }\right) \right) $ is continuous. Moreover for
all $h>0$ with $h\leq $ $h_{0}<\infty $ and $u\in W_{p,\gamma }^{\left[ m%
\right] }\left( \Omega ;E\left( A\right) ,E\right) $ the following estimate
holds 
\begin{equation*}
\left\Vert D^{\left[ \alpha \right] }u\right\Vert _{L_{p}\left( \Omega
;E\left( A^{1-\varkappa -\mu }\right) \right) }\leq h^{\mu }\left\Vert
u\right\Vert _{W_{p,\gamma }^{\left[ m\right] }\left( \Omega ;E\left(
A\right) ,E\right) }+h^{-\left( 1-\mu \right) }\left\Vert u\right\Vert
_{L_{p}\left( \Omega ;E\right) }.
\end{equation*}

\begin{center}
\bigskip \textbf{2. Singular degenerate elliptic DOE}
\end{center}

\ Consider the BVP for the following singular degenerate DOE 
\begin{equation}
-\dsum\limits_{k=1}^{n}\left[ x_{k}^{2\alpha _{k}}\frac{\partial ^{2}u}{%
\partial x_{k}^{2}}+x_{k}^{\alpha _{k}}A_{k}\left( x\right) \frac{\partial u%
}{\partial x_{k}}\right] +Au+\lambda u=f\left( x\right) ,\text{ }x\in G 
\tag{2.1}
\end{equation}

\begin{equation}
L_{k}u=\sum\limits_{i=0}^{m_{k}}\left[ \delta _{ki}u_{x_{k}}^{\left[ i\right]
}\left( a_{k},x\left( k\right) \right) +\sum\limits_{j=0}^{N_{k}}\nu
_{kij}u_{x_{k}}^{\left[ i\right] }\left( x_{kij},x\left( k\right) \right) %
\right] =0\text{, }  \tag{2.2}
\end{equation}%
where $x\left( k\right) \in G_{k}$ and%
\begin{eqnarray*}
u_{x_{k}}^{\left[ i\right] } &=&\left[ x_{k}^{\alpha _{k}}\frac{\partial }{%
\partial x_{k}}\right] ^{i}u\left( x\right) ,\text{ }G=\dprod%
\limits_{k=1}^{n}\left( 0,a_{k}\right) ,\text{ }G_{k}=\dprod\limits_{j\neq
k}\left( 0,a_{j}\right) ,\text{ }\delta _{km_{k}}\neq 0 \\
m_{k} &\in &\left\{ 0,1\right\} \text{, }x\left( k\right) =\left(
x_{1},x_{2},...,x_{k-1},x_{k+1},...,x_{n}\right) \text{, }j,\text{ }%
k=1,2,...,n;
\end{eqnarray*}%
$\delta _{ki},$ $\nu _{kij}$ are complex numbers, $\lambda $ is a complex
parameter, $A$ and $A_{k}\left( x\right) $ are linear operators in a Banach
space $E.$

Let $\alpha =\left( \alpha _{1},\alpha _{2},...,\alpha _{n}\right) ,$ $%
\gamma _{k}\left( x\right) =x_{k}^{\alpha _{k}}.$ The main result is the
following

\textbf{Theorem 2.1. }Assume the following conditions are satisfied:

(1) $E$ is an UMD space and $A$ is a $R$-positive operator in $E;$

(2) $1+\frac{1}{p}<\alpha _{k}<\frac{p-1}{2},$ $p\in \left( 1,\infty \right)
,$ $\delta _{km_{k}}\neq 0$ $;$

(3) for any $\varepsilon >0$, there is a positive $C\left( \varepsilon
\right) $ such that

$\left\Vert A_{k}\left( y\right) u\right\Vert \leq \varepsilon \left\Vert
u\right\Vert _{\left( E\left( A\right) ,E\right) _{\frac{1}{2},\infty
}}+C\left( \varepsilon \right) \left\Vert u\right\Vert $ for $u\in \left(
E\left( A\right) ,E\right) _{\frac{1}{2},\infty }.$

\vspace{3mm}Then, the problem $\left( 2.1\right) -\left( 2.2\right) $ has a
unique solution $u\in W_{p,\alpha }^{2}\left( G;E\left( A\right) ,E\right) $
for $f\in L_{p}\left( G;E\right) $ and sufficiently large $\left\vert
\lambda \right\vert $ with $\left\vert \arg \lambda \right\vert \leq \varphi 
$ and the following uniform coercive estimate holds

\begin{equation}
\sum\limits_{k=1}^{n}\sum\limits_{i=0}^{2}\left\vert \lambda \right\vert ^{1-%
\frac{i}{2}}\left\Vert x_{k}^{i\alpha _{k}}\frac{\partial ^{i}u}{\partial
x_{k}^{i}}\right\Vert _{L_{p}\left( G;E\right) }+\left\Vert Au\right\Vert
_{L_{p}\left( G;E\right) }\leq M\left\Vert f\right\Vert _{L_{p}\left(
G;E\right) }.  \tag{2.3}
\end{equation}

For proving the main theorem, consider at first the BVP for the singular
degenerate ordinary DOE

\begin{equation}
\ -u^{\left[ 2\right] }\left( t\right) +\left( A+\lambda \right) u\left(
t\right) =f,\text{ }t\in \left( 0,a\right) ,  \tag{2.4}
\end{equation}%
\begin{equation*}
L_{1}u=\sum\limits_{i=0}^{m}\left[ \delta _{i}u^{\left[ i\right] }\left(
a\right) +\sum\limits_{j=1}^{N}\nu _{ij}u^{\left[ i\right] }\left(
t_{ij}\right) \right] =0,
\end{equation*}%
where $u^{\left[ i\right] }=\left( t^{\nu }\frac{d}{dt}\right) ^{i}$, $\nu
>1,$ $m\in \left\{ 0,1\right\} ;$ $\delta _{i},$ $\nu _{ij},$ are complex
numbers and $x_{ij}\in \left( 0,a\right) ;$ $A$ is a possible unbounded
operators in $E.$

\textbf{Remark 2.1. }Let 
\begin{equation}
\tau =-\int\limits_{t}^{a}z^{-\nu }dz,\text{ }t=\left[ a^{1-\nu }-\left( \nu
-1\tau \right) \right] ^{\frac{1}{\nu -1}}.  \tag{2.5}
\end{equation}

Under the substitution $\left( 2.5\right) $ the spaces $L_{p}\left(
0,a;E\right) $, $W_{p,\nu }^{\left[ 2\right] }\left( 0,a;E\left( A\right)
,E\right) $ are mapped isomorphically onto weighted spaces%
\begin{equation*}
L_{p,\tilde{\nu}}\left( -\infty ,0;E\right) ,W_{p,\tilde{\nu}}^{2}\left(
-\infty ,0;E\left( A\right) ,E\right) ,
\end{equation*}
respectively, where $\tilde{\nu}=\nu \left( t\left( \tau \right) \right) .$
Moreover, under the substitution $\left( 2.5\right) $ the problem $\left(
2.1\right) -\left( 2.2\right) $ is transformed into the following non
degenerate problem

\begin{equation}
-u^{\left( 2\right) }\left( \tau \right) +Au\left( \tau \right) =f\left(
\tau \right) ,  \tag{2.6}
\end{equation}%
\begin{equation*}
L_{1}u=\sum\limits_{i=0}^{m}\left[ \delta _{i}u^{\left( i\right) }\left(
a\right) +\sum\limits_{j=1}^{N}\nu _{ij}u^{\left( i\right) }\left( \tau
_{ij}\right) \right] =0
\end{equation*}%
considered in the weighted space $L_{p,\tilde{\nu}}\left( -\infty
;0;E\right) .$

In a similar way as in $\left[ \text{13, Theorem 4.1}\right] $ and $\left[ 
\text{11, Lemma 3.2}\right] ,$ we obtain

\textbf{Proposition 2.}1\textbf{. }Let the following conditions be satisfied:

(1)\ $E$ is a UMD\ space Banach space and $A$ is an $R$ positive in $E:$

(2) $1+\frac{1}{p}<\nu <\frac{\left( p-1\right) }{2},1<p<\infty ,$ $\delta
_{m}\neq 0.$

Then, the problem $\left( 2.4\right) $ has a unique solution $u\in W_{p,\nu
}^{\left[ 2\right] }\left( 0,a;E\left( A\right) ,E\right) $ for all \ $f\in
L_{p}\left( 0,a;E\right) $, for $\left\vert \arg \lambda \right\vert \leq
\varphi $ with sufficiently large $\left\vert \lambda \right\vert $ and the
uniform coercive estimate holds

\begin{equation*}
\sum\limits_{i=0}^{2}\left\vert \lambda \right\vert ^{1-\frac{i}{2}%
}\left\Vert u^{\left[ i\right] }\right\Vert _{L_{p}\left( 0,a;E\right)
}+\left\Vert Au\right\Vert _{L_{p}\left( 0,a;E\right) }\leq C\left\Vert
f\right\Vert _{L_{p}\left( 0,a;E\right) }.
\end{equation*}

\textbf{Proof. }Consider the transformed problem $\left( 2.6\right) $. Since
the operator $A$ generates an analytic semigroups, by reasoning as in $\left[
\text{4, Lemma 5. 3. 2/1}\right] $ we fined the representation of solution
of this problem. Then by using the properties of positive operator $A,$ the
estimates of analytic semigroups and integral operators in weighted spaces $%
L_{p,\tilde{\nu}}\left( -\infty ,0;E\right) $ we obtain the assertion.

Consider the operator $B$ generated by problem $\left( 2.4\right) $, i.e. 
\begin{equation*}
D\left( B\right) =W_{p,\nu }^{\left[ 2\right] }\left( 0,a;E\left( A\right)
,E,L_{1}\right) \text{, }Bu=-u^{\left[ 2\right] }+Au.
\end{equation*}

In a similar way as in $\left[ \text{11, Theorem 3.1}\right] $ we obtain

\textbf{Proposition 2.2}$.$ Suppose all conditions of Proposition 2.1 are
satisfied. Then, the operator $B$ is $R$-positive in $L_{p}\left(
0,a;E\right) .$

Proposition 2.1 implies that the operator $B$ is positive in $L_{p}\left(
0,a;E\right) $ and also is a generator of an analytic semigroup. Consider
the principal part of the problem $\left( 2.4\right) ,$ i.e. consider the
problem $\left( 2.6\right) $.

\textbf{Proposition 2.3}$.$ Assume all conditions of the Proposition 2.1 are
satisfied. Then, the problem $\left( 2.6\right) $ has a unique solution $%
u\in W_{p,\nu }^{2}\left( 0,a;E\left( A\right) ,E\right) $ for all \ $f\in
L_{p}\left( 0,a;E\right) $, $\left\vert \arg \lambda \right\vert \leq
\varphi $ and sufficiently large $\left\vert \lambda \right\vert .$
Moreover, the uniform coercive estimate holds

\begin{equation}
\sum\limits_{i=0}^{2}\left\vert \lambda \right\vert ^{1-\frac{i}{2}%
}\left\Vert \tau ^{i\nu }u^{\left( i\right) }\right\Vert _{L_{p}\left(
0,a;E\right) }+\left\Vert Au\right\Vert _{L_{p}\left( 0,a;E\right) }\leq
C\left\Vert f\right\Vert _{L_{p}\left( 0,a;E\right) }.  \tag{2.7}
\end{equation}

\textbf{Proof. \ }Since $\nu >1,$ by Theorem B and Remark 2.1 for all $%
\varepsilon >0$ there is a continuous function $C\left( \varepsilon \right) $
such that 
\begin{equation}
\left\Vert \nu x^{\nu -1}u^{\left[ 1\right] }\right\Vert _{L_{p}\left(
0,a;E\right) }\leq \varepsilon \left\Vert u\right\Vert _{W_{p,\nu }^{\left[ 2%
\right] }\left( 0,a;E\left( A\right) ,E\right) }+C\left( \varepsilon \right)
\left\Vert u\right\Vert _{L_{p}\left( 0,a;E\right) }.  \tag{2.8}
\end{equation}%
Then, in view of $\left( 2.6\right) ,\left( 2.7\right) $ and due to
positivity of operator $B$ we have the following estimate 
\begin{equation}
\left\Vert \nu x^{\nu -1}u^{\left[ 1\right] }\right\Vert _{L_{p}\left(
0,a;E\right) }\leq \varepsilon \left\Vert Bu\right\Vert _{L_{p}\left(
0,a;E\right) }.  \tag{2.9}
\end{equation}

S\i nce $-x^{2\nu }u^{\left( 2\right) }=-u^{\left[ 2\right] }+\nu x^{\nu
-1}u^{\left[ 1\right] },$ the assertion is obtained from Proposition 2.1 and
estimate $\left( 2.9\right) .$

Consider the operator $S$ generated by problem $\left( 2.6\right) $, i.e. 
\begin{equation*}
D\left( S\right) =W_{p,\nu }^{2}\left( 0,a;E\left( A\right) ,E,L_{k}\right) ,%
\text{ }Su=-x^{2\nu }u^{\left( 2\right) }+Au.
\end{equation*}

\ \textbf{Result 2.1}$.$Suppose all conditions of Proposition 2.1 are
satisfied. Then, the operator $S$ is $R$-positive in $L_{p}\left(
0,a;E\right) .$

\textbf{\ }The assertion is obtained from Proposition 2.2 and the estimate $%
\left( 2.9\right) .$

Consider now the principal part of the problem $\left( 2.1\right) -\left(
2.2\right) $ with constant coefficients, i.e. 
\begin{equation}
-\dsum\limits_{k=1}^{n}x_{k}^{2\alpha _{k}}\frac{\partial ^{2}u}{\partial
x_{k}^{2}}+Au+\lambda u=f\left( x\right) ,\text{ }L_{k}u=0\text{, }%
k=1,2,...,n  \tag{2.10}
\end{equation}

\textbf{Proposition 2.4. }Assume $E$ is a UMD space and $A$ is an $R$%
-positive operator in $E.$\ Let $1+\frac{1}{p}<\alpha _{k}<\frac{p-1}{2},$ $%
p\in \left( 1,\infty \right) $.

\vspace{3mm}Then problem $\left( 2.10\right) $ has a unique solution $u\in
W_{p,\alpha }^{2}\left( G;E\left( A\right) ,E\right) $ for $f\in L_{p}\left(
G;E\right) $ and sufficiently large $\left\vert \lambda \right\vert $ with $%
\left\vert \arg \lambda \right\vert \leq \varphi $ and the uniform coercive
estimate holds

\begin{equation}
\sum\limits_{k=1}^{n}\sum\limits_{i=0}^{2}\left\vert \lambda \right\vert ^{1-%
\frac{i}{2}}\left\Vert x_{k}^{i\alpha }\frac{\partial ^{i}u}{\partial
x_{k}^{i}}\right\Vert _{L_{p}\left( G;E\right) }+\left\Vert Au\right\Vert
_{L_{p}\left( G;E\right) }\leq M\left\Vert f\right\Vert _{L_{p}\left(
G;E\right) }.  \tag{2.11}
\end{equation}

\textbf{Proof. }Consider first all of the problem $\left( 2.10\right) $ for $%
n=2$ i.e

\begin{equation}
-\dsum\limits_{k=1}^{2}x_{k}^{2\alpha _{k}}\frac{\partial ^{2}u}{\partial
x_{k}^{2}}+Au+\lambda u=f\left( x_{1},x_{2}\right) ,\text{ }L_{k}u=0\text{, }%
k=1,2.  \tag{2.12}
\end{equation}%
Since%
\begin{equation*}
L_{p}\left( 0,a_{2};L_{p}\left( 0,a_{1};E\right) \right) =L_{p}\left( \left(
0,a_{1}\right) \left( 0,a_{2}\right) \times ;E\right)
\end{equation*}%
then the BVP $\left( 2.12\right) $ can be expressed as

\begin{equation}
-x^{2\alpha _{2}}\frac{d^{2}u}{dx_{2}^{2}}+\left( S+\lambda \right) u\left(
x_{2}\right) =f\left( x_{2}\right) \text{, }L_{2}u=0.  \tag{2.13}
\end{equation}

By virtue of $\left[ \text{1, Theorem 4.5.2}\right] $, $F=L_{p}\left(
0,a_{1};E\right) \in UMD$ provided $E\in UMD$, $p\in \left( 1,\infty \right) 
$. Then, by virtue of $\left[ 21\right] $, $L_{p}\left( 0,a_{1};E\right) $
is the space satisfying the multiplier condition. By Result 2.1\ the
operator $S\ $is $R$-positive in $F.$ Then by virtue of Proposition 2.3 we
get that, for $f\in L_{p}\left( 0,a_{2};F\right) $ the problem $\left(
2.13\right) ,$ i.e. problem $\left( 2.12\right) $ for $\left\vert \arg
\lambda \right\vert \leq \varphi $ and sufficiently large $\left\vert
\lambda \right\vert $ has a unique solution $u\in $ $W_{p,\alpha
_{2}}^{2}\left( 0,a_{2};D\left( S\right) ,F\right) $ and the coercive
uniform estimate $\left( 2.11\right) $ holds for solution of the problem $%
\left( 2.12\right) $. By continuing the above proses $n$ time, we obtain
that the problem $\left( 2.10\right) $\ has a unique solution $u\in
W_{p,\alpha }^{2}\left( G;E\left( A\right) ,E\right) $ for $f\in L_{p}\left(
G;E\right) $ and $\left\vert \arg \lambda \right\vert \leq \varphi $ and
sufficiently large $\left\vert \lambda \right\vert ,$ moreover, the uniform
estimate $\left( 2.11\right) $ holds.

\textbf{Proof of Theorem 2.1.} Let $Q_{0}$ denote differential operator
generated by problem $\left( 9\right) $ i.e., 
\begin{equation*}
D\left( Q_{0}\right) =\left\{ u\in W_{p,\alpha }^{2}\left( G;E\left(
A\right) ,E\right) ,\text{ }L_{k}u=0,\text{ }k=1,2,...,n\right\} ,
\end{equation*}%
\begin{equation*}
Q_{0}u=-\dsum\limits_{k=1}^{2}x_{k}^{2\alpha _{k}}\frac{\partial ^{2}u}{%
\partial x_{k}^{2}}+Au.
\end{equation*}%
The estimate of type $\left( 2.11\right) $ implies that the operator $Q_{0}$
has a bounded inverse from $L_{p}\left( G;E\right) $ to $W_{p,\alpha
}^{2}\left( G;E\left( A\right) ,E\right) $, i.e. for all $f\in L_{p}\left(
G;E\right) ,$ $\lambda \in S\left( \varphi \right) $ with sufficiently large 
$\left\vert \lambda \right\vert $ the estimate holds%
\begin{equation*}
\left\Vert \left( Q_{0}+\lambda \right) ^{-1}f\right\Vert _{W_{p,\alpha
}^{2}\left( G;E\left( A\right) ,E\right) }\leq C\left\Vert f\right\Vert
_{L_{p}\left( G;E\right) }.
\end{equation*}%
Moreover, by virtue of Theorem B and in view of assumption (3), for all $%
\varepsilon >0$ there is a continuous function $C\left( \varepsilon \right) $
such that 
\begin{equation*}
\dsum\limits_{k=1}^{n}\left\Vert x_{k}^{\alpha _{k}}A_{k}u\right\Vert
_{L_{p}\left( G;E\right) }\leq \varepsilon \left\Vert u\right\Vert
_{W_{p,\alpha }^{2}\left( G;E\left( A\right) ,E\right) }+C\left( \varepsilon
\right) \left\Vert u\right\Vert _{L_{p}\left( G;E\right) }.
\end{equation*}%
From the above estimates we obtain that there is a positive number $\delta
<1 $ such that 
\begin{equation*}
\left\Vert Q_{1}u\right\Vert _{L_{p}\left( G;E\right) }<\delta \left\Vert
\left( Q_{0}+\lambda \right) u\right\Vert _{L_{p}\left( G;E\right) }
\end{equation*}%
for all $u\in W_{p,\alpha }^{2}\left( G;E\left( A\right) ,E\right) ,$ where 
\begin{equation*}
Q_{1}u=\dsum\limits_{k=1}^{n}x_{k}^{\alpha _{k}}A_{k}\left( x\right) \frac{%
\partial u}{\partial x_{k}}.
\end{equation*}

Let $Q$ denote differential operator generated by problem $\left( 2.1\right)
-\left( 2.2\right) $ for $\lambda =0.$ It is clear that 
\begin{equation*}
\left( Q+\lambda \right) =\left[ I+Q_{1}\left( Q_{0}+\lambda \right) ^{-1}%
\right] \left( Q_{0}+\lambda \right) .
\end{equation*}%
Therefore, we obtain that the operator $\left( Q+\lambda \right) ^{-1}$ is
bounded from $L_{p}\left( G;E\right) $ to $W_{p,\alpha }^{2}\left( G;E\left(
A\right) ,E\right) $ and the estimate $\left( 2.11\right) $ is satisfied.

Let $B=L\left( L_{p}\left( G;E\right) \right) .$ We get the following result
from Theorem 2.1:

\textbf{Result 2.2. }Theorem 2.1 implies that differential operator $Q$ has
a resolvent $\left( Q+\lambda \right) ^{-1}$ for $\left\vert \arg \lambda
\right\vert \leq \varphi $ and the following estimate holds 
\begin{equation*}
\sum\limits_{k=1}^{n}\sum\limits_{i=0}^{2}\left\vert \lambda \right\vert ^{1-%
\frac{i}{2}}\left\Vert x_{k}^{i\alpha }\frac{\partial ^{i}}{\partial
x_{k}^{i}}\left( Q+\lambda \right) ^{-1}\right\Vert _{B}+\left\Vert A\left(
Q+\lambda \right) ^{-1}\right\Vert _{B}\leq M.
\end{equation*}

\textbf{Proposition 2.5}$.$ Assume\ $E$ is a UMD space and $A$ is a $R$
positive in $E$. Let $1+\frac{1}{p}<\alpha _{k}<\frac{p-1}{2},1<p<\infty ,$ $%
\delta _{m}\neq 0$. Then, the operator $Q_{0}$ is $R$-positive in $%
L_{p}\left( G;E\right) .$

\textbf{Proof. }By\ reasoning as in the proof of Theorem 2.1, we get that
problem $\left( 2.10\right) $ can be expressed as the problem $\left(
2.13\right) .$ By virtue of Result 2.1 the operator $S$ is $R$-positive,
then by applying again the Result 2.1 to problem $\left( 2.13\right) $ in $%
L_{p}\left( 0,a_{2};F\right) ,$ $F=L_{p}\left( 0,a_{1};E\right) $ and by
continuing it $n$ time we obtain that the operator $Q_{0}$ is $R$-positive
in $L_{p}\left( G;E\right) .$

\textbf{Remark 2.2. }Note that, by using the techniques similar to those
applied in Theorem 2.1 we obtain the same results for differential-operator
equations of the arbitrary order.

\begin{center}
\textbf{3. Cauchy problem for singular degenerate parabolic equation }
\end{center}

Consider the mixed problem for singular degenerate parabolic DOE equation 
\begin{equation}
\frac{\partial u}{\partial t}-\dsum\limits_{k=1}^{n}x_{k}^{2\alpha _{k}}%
\frac{\partial ^{2}u}{\partial x_{k}^{2}}+Au=f\left( t,x\right) ,\text{ }%
x\in G,\text{ }t\in R_{+},  \tag{3.1}
\end{equation}%
\begin{equation}
\sum\limits_{i=0}^{m_{k}}\delta _{ki}u_{x}^{\left[ i\right] }\left(
t,a_{k},x\left( k\right) \right) =0,\text{ }x\left( k\right) \in G_{k}, 
\tag{3.2}
\end{equation}

\begin{equation*}
\text{  }u\left( 0,x\right) =0\text{, }
\end{equation*}%
where $\delta _{ki}$ are complex numbers, $G=\dprod\limits_{k=1}^{n}\left(
0,a_{k}\right) ,$ $G_{k}=\dprod\limits_{j\neq k}\left( 0,a_{j}\right) $,  $%
j=1,2,...,n,$ $x\left( k\right) =\left(
x_{1},x_{2},...,x_{k-1},x_{k+1},...,x_{n}\right) ,$ $A$ is a possible
unbounded operator in a Banach space $E$ and $u=u\left( t,x\right) $ is a
solution of the problem $\left( 3.1\right) -\left( 3.2\right) .$

For $\mathbf{p=}\left( p,p_{1}\right) $, $\Delta _{+}=R_{+}\times G,$ $L_{%
\mathbf{p}}\left( \Delta _{+};E\right) $ will be denoted the space of all $E$%
-valued $\mathbf{p}$-summable functions with mixed norm   i.e. the space of
all measurable functions $f$ defined on $\Delta _{+}$, for which 
\begin{equation*}
\left\Vert f\right\Vert _{L_{\mathbf{p,\alpha }}\left( \Delta _{+}\right)
}=\left( \int\limits_{R_{+}}\left( \dint\limits_{G}\left\Vert f\left(
x\right) \right\Vert ^{p}\alpha \left( x\right) dx\right) ^{\frac{p_{1}}{p}%
}dt\right) ^{\frac{1}{p_{1}}}<\infty .
\end{equation*}%
Analogously, $W_{\mathbf{p,}\alpha }^{m}\left( \Delta _{+},E\left( A\right)
,E\right) $ denotes the Sobolev space with corresponding mixed norm (see $%
\left[ \text{23, \S }\right] $ for scalar case)$.$

\textbf{Theorem 3.1. }Assume all conditions of Theorem 2.1 hold for $\varphi
>\frac{\pi }{2}$. Then, for all $f\in L_{\mathbf{p}}\left( \Delta
_{+};E\right) $ and sufficiently large $d>0$\ problem $\left( 11\right)
-\left( 12\right) $ has a unique solution belonging to $W_{\mathbf{p},\alpha
}^{1,2}\left( \Delta _{+};E\left( A\right) ,E\right) $ and the following
coercive estimate holds 
\begin{equation*}
\left\Vert \frac{\partial u}{\partial t}\right\Vert _{L_{\mathbf{p}}\left(
G_{+};E\right) }+\dsum\limits_{k=1}^{n}\left\Vert x^{2\alpha _{k}}\frac{%
\partial ^{2}u}{\partial x_{k}^{2}}\right\Vert _{L_{\mathbf{p}}\left(
G_{+};E\right) }+\left\Vert Au\right\Vert _{L_{\mathbf{p}}\left(
G_{+};E\right) }\leq
\end{equation*}%
\begin{equation*}
C\left\Vert f\right\Vert _{L_{\mathbf{p}}\left( G_{+};E\right) }.
\end{equation*}%
\textbf{Proof. }The problem $\left( 3.1\right) -\left( 3.2\right) $ can be
expressed as the following Cauchy problem 
\begin{equation}
\frac{du}{dt}+Q_{0}u\left( t\right) =f\left( t\right) ,\text{ }u\left(
0\right) =0.  \tag{3.3}
\end{equation}%
Preposition 2.4 implies that the operator $Q_{0}$ is $R$-positive in $%
F=L_{p}\left( G;E\right) .$ By $\left[ \text{22, \S 1.14}\right] ,$ $Q_{0}$
is a generator of an analytic semigroup in $F.$ Then, by virtue of $\left[ 
\text{16, Theorem 4.2}\right] ,$ we obtain that for all $f\in
L_{p_{1}}\left( R_{+};F\right) $ and sufficiently large $d>0,$ problem $%
\left( 15\right) $ has a unique solution belonging to $W_{p_{1}}^{1}\left(
R_{+};D\left( Q_{0}\right) ,F\right) $ and the estimate holds 
\begin{equation*}
\left\Vert \frac{du}{dt}\right\Vert _{L_{p_{1}}\left( R_{+};F\right)
}+\left\Vert Qu\right\Vert _{L_{p_{1}}\left( R_{+};F\right) }\leq
C\left\Vert f\right\Vert _{L_{p_{1}}\left( R_{+};F\right) }.
\end{equation*}

Since $L_{p_{1}}\left( R_{+};F\right) =L_{\mathbf{p}}\left( \Delta
_{+};E\right) ,$ by Theorem 2.1, we have 
\begin{equation*}
\left\Vert Q_{0}u\right\Vert _{L_{p_{1}}\left( R_{+};F\right) }=D\left(
Q_{0}\right) .
\end{equation*}%
These relations and the above estimate prove the hypothesis to be true.

\begin{center}
\textbf{\ 4. Singular degenerate boundary value problems for infinite
systems of equations }
\end{center}

Consider the infinite system of BVPs%
\begin{equation}
-x^{2\alpha }\frac{\partial ^{2}u_{m}}{\partial x^{2}}-x^{2\beta }\frac{%
\partial ^{2}u_{m}}{\partial y^{2}}+d_{m}u_{m}+\dsum\limits_{j=1}^{\infty
}x^{\alpha }a_{mj}\left( x,y\right) \frac{\partial u_{j}}{\partial x} 
\tag{4.1}
\end{equation}%
\begin{equation*}
+\dsum\limits_{j=1}^{\infty }y^{\beta }b_{mj}\left( x,y\right) \frac{%
\partial u_{j}}{\partial y}+\lambda u=f_{m}\left( x,y\right) ,\text{ }%
L_{1}u=0,\text{ }L_{2}u=0,
\end{equation*}%
where $L_{k}$ are defined by $\left( 2.2\right) $ and $x,y\in G=\left(
0,a\right) \times \left( 0,b\right) $. 
\begin{equation*}
D=\left\{ d_{m}\right\} ,\text{ }d_{m}>0,\text{ }u=\left\{ u_{m}\right\} ,%
\text{ }Du=\left\{ d_{m}u_{m}\right\} ,\text{ }m=1,2,...,
\end{equation*}

\begin{equation*}
\text{ }l_{q}\left( D\right) =\left\{ u\text{: }u\in l_{q},\right.
=\left\Vert u\right\Vert _{l_{q}\left( D\right) }=\left. \left(
\sum\limits_{m=1}^{\infty }\left\vert d_{m}u_{m}\right\vert ^{q}\right) ^{%
\frac{1}{q}}<\infty ,q\in \left( 1,\infty \right) \right\} .
\end{equation*}%
From Theorem 2.1, we obtain

\textbf{Theorem 4.1. }Assume$\ p\in \left( 1,\infty \right) ,$ $1<\alpha
,\beta <p-1,$ $a_{mj},b_{mj}\in L_{\infty }\left( G\right) $. For $0<\mu <%
\frac{1}{2},$ $0<\nu <1$\ and for all $x,y\in G$ 
\begin{equation*}
\text{ }\sup\limits_{m}\sum\limits_{j=1}^{\infty }a_{mj}\left( x,y\right)
d_{j}^{-\left( \frac{1}{2}-\mu \right) }<M,\text{ }\sup\limits_{m}\sum%
\limits_{j=1}^{\infty }b_{mj}\left( x,y\right) d_{j}^{-\left( 1-\nu \right)
}.
\end{equation*}%
Then, for all $f\left( x,y\right) =\left\{ f_{m}\left( x,y\right) \right\}
_{1}^{\infty }\in L_{p}\left( \left( G\right) ;l_{q}\right) ,$ $p,q\in
\left( 1,\infty \right) $, $\left\vert \arg \lambda \right\vert \leq \varphi 
$, $0\leq \varphi <\pi $ and for sufficiently large $\left\vert \lambda
\right\vert ,$ problem $\left( 4.1\right) $ has a unique solution $u=\left\{
u_{m}\left( x,y\right) \right\} _{1}^{\infty }$ that belongs to space $%
W_{p,\alpha ,\beta }^{2}\left( G,l_{q}\left( D\right) ,l_{q}\right) $ and 
\begin{equation*}
\sum\limits_{i=0}^{2}\left\vert \lambda \right\vert ^{1-\frac{i}{2}}\left[
\left\Vert x^{i\alpha }\frac{\partial ^{i}u}{\partial x^{i}}\right\Vert
_{L_{p}\left( G;l_{q}\right) }+\left\Vert y^{i\beta }\frac{\partial ^{i}u}{%
\partial y^{i}}\right\Vert _{L_{p}\left( G;l_{q}\right) }\right] 
\end{equation*}

\begin{center}
\ 
\begin{equation*}
+\left\Vert Du\right\Vert _{L_{p}\left( G;l_{q}\right) }\leq M\left\Vert
f\right\Vert _{L_{p}\left( G;l_{q}\right) }.
\end{equation*}

\textbf{5. Wentzell-Robin type mixed problem for degenerate parabolic
equation }
\end{center}

Consider the problem%
\begin{equation}
\frac{\partial u}{\partial t}-\dsum\limits_{k=1}^{n}x_{k}^{2\alpha _{k}}%
\frac{\partial ^{2}u}{\partial x_{k}^{2}}+a\frac{\partial ^{2}u}{\partial
y^{2}}+b\frac{\partial u}{\partial y}+cu=f\left( x,y,t\right) ,\text{ } 
\tag{5.1}
\end{equation}%
\ \ \ 

\begin{equation}
L_{k}u=0\text{, }B_{j}u=0\text{, }j=0,1,\text{ }t\in \left( 0,T\right) \text{%
, }x\in \sigma ,\text{ }y\in \left( 0,1\right) ,  \tag{5.2}
\end{equation}

\begin{equation}
u\left( x,y,0\right) =0\text{, }x\in G,\text{ }y\in \left( 0,1\right) ,\text{
}  \tag{5.3}
\end{equation}%
where $a=a\left( x,y,t\right) ,$ $a_{1}=a_{1}\left( x,y,t\right) ,$ $%
b_{1}=b_{1}\left( x,y,t\right) ,$ $c=c\left( x,y,t\right) $ are
complex-valued functions on $\tilde{\Omega}=G\times \left( 0,1\right) \times
\left( 0,T\right) $. For $\mathbf{\tilde{p}=}\left( p,p_{1},2\right) $ and $%
L_{\mathbf{\tilde{p}}}\left( \tilde{\Omega}\right) $ will denote the space
of all $\mathbf{\tilde{p}}$-summable scalar-valued\ functions with mixed
norm. Analogously, $W_{\mathbf{\tilde{p},}\alpha ,\beta }^{2,1}\left( \tilde{%
\Omega}\right) $ denotes the Sobolev space with corresponding mixed norm,
i.e., $W_{\mathbf{\tilde{p},}\alpha ,\beta }^{2,1}\left( \tilde{\Omega}%
\right) $ denotes the space of all functions $u\in L_{\mathbf{\tilde{p}}%
}\left( \tilde{\Omega}\right) $ possessing the derivatives $\frac{\partial u%
}{\partial t},$ $x^{2\alpha }\frac{\partial ^{2}u}{\partial x^{2}},$ $%
y^{2\beta }\frac{\partial ^{2}u}{\partial y^{2}}\in L_{\mathbf{\tilde{p}}%
}\left( \tilde{\Omega}\right) $ with the norm 
\begin{equation*}
\ \left\Vert u\right\Vert _{W_{\mathbf{\tilde{p},}\alpha ,\beta
}^{2,1}\left( \tilde{\Omega}\right) }=\left\Vert u\right\Vert _{L_{\mathbf{%
\tilde{p}}}\left( \tilde{\Omega}\right) }+\left\Vert \frac{\partial u}{%
\partial t}\right\Vert _{L_{\mathbf{\tilde{p}}}\left( \tilde{\Omega}\right)
}+\left\Vert x^{2\alpha }\frac{\partial ^{2}u}{\partial x^{2}}\right\Vert
_{L_{\mathbf{\tilde{p}}}\left( \tilde{\Omega}\right) }+\left\Vert y^{2\beta }%
\frac{\partial ^{2}u}{\partial y^{2}}\right\Vert _{L_{\mathbf{\tilde{p}}%
}\left( \tilde{\Omega}\right) }.
\end{equation*}

\textbf{Condition 5.1 }Assume;

(1) $1<\alpha ,$ $\beta <p-1,$ $p,$ $p_{1}\in \left( 1,\infty \right) ;$

(2)\ $a_{1}\left( x,.,t\right) \in W_{\infty }^{1}\left( 0,1\right) ,$ $%
a_{1}\left( x,.,t\right) \geq \delta >0,$ $b_{1}\left( x,.,t\right) ,$ $%
c\left( x,.,t\right) \in L_{\infty }\left( 0,1\right) $ for a.e. $x\in G,$ $%
t\in \left( 0,T\right) ;$

(3) $b\left( .,y,t\right) ,$ $c\left( .,y,t\right) \in C\left( \bar{G}%
\right) $ for $y\in \left( 0,b\right) $ and $t\in \left( 0,T\right) ;$

\bigskip\ In this section, we present the following result:

\bigskip \bigskip \textbf{Theorem 5.1. }Suppose the Condition 5.1 hold.
Then, for $f\in L_{\mathbf{\tilde{p}}}\left( \tilde{\Omega};E\right) $
problem $\left( 5.1\right) -\left( 5.3\right) $ has a unique solution $u$

belonging to $W_{\mathbf{\bar{p},}\alpha ,\beta }^{2,1}\left( \tilde{\Omega}%
\right) $ and the following coercive estimate holds 
\begin{equation*}
\left\Vert \frac{\partial u}{\partial t}\right\Vert _{L_{\mathbf{\bar{p}}%
}\left( \tilde{\Omega};E\right) }+\left\Vert x^{2\alpha }\frac{\partial ^{2}u%
}{\partial x^{2}}\right\Vert _{L_{\mathbf{\tilde{p}}}\left( \tilde{\Omega}%
\right) }+\left\Vert y^{2\beta }\frac{\partial ^{2}u}{\partial y^{2}}%
\right\Vert _{L_{\mathbf{\tilde{p}}}\left( \tilde{\Omega}\right)
}+\left\Vert Au\right\Vert _{L_{\mathbf{\bar{p}}}\left( G_{T};E\right) }\leq
C\left\Vert f\right\Vert _{L_{\mathbf{\bar{p}}}\left( \tilde{\Omega}%
;E\right) }.
\end{equation*}

\ \textbf{Proof.} Let $E=L_{2}\left( 0,1\right) $. It is known $\left[ 10%
\right] $\ that $L_{2}\left( 0,1\right) $ is an $UMD$ space. Consider the
operator $A$ defined by 
\begin{equation*}
D\left( A\right) =W_{2}^{2}\left( 0,1;B_{j}u=0\right) ,\text{ }Au=a_{1}\frac{%
\partial ^{2}u}{\partial y^{2}}+b_{1}\frac{\partial u}{\partial y}+cu.
\end{equation*}

Therefore, the problem $\left( 5.1\right) -\left( 5.3\right) $ can be
rewritten in the form of $\left( 3.1\right) -\left( 3.2\right) $, where $%
u\left( x\right) =u\left( x,.\right) ,$ $f\left( x\right) =f\left(
x,.\right) $\ are functions with values in $E=L_{2}\left( 0,1\right) .$ By
virtue of $\left[ \text{24, 25}\right] $ the operator $A$ generates analytic
semigroup in $L_{2}\left( 0,b\right) $. Then in view of Hill-Yosida theorem
(see e.g. $\left[ \text{22, \S\ 1.13}\right] $) this operator is positive in 
$L_{2}\left( 0,b\right) .$ Since all uniform bounded set in Hilbert space is 
$R$-bounded (see $\left[ 4\right] $ ), i.e. we get that the operator $A$ is $%
R$-positive in $L_{2}\left( 0,b\right) .$ Then from Theorem 3.1 we obtain
the assertion.

\bigskip \textbf{References}

\begin{quote}
\ \ \ \ \ \ \ \ \ \ \ \ \ \ \ \ \ \ \ \ \ \ \ \ \ \ \ \ \ \ \ \ \ \ \ \ \ \
\ \ \ \ \ \ \ \ \ \ \ \ \ \ \ \ \ \ \ \ \ \ \ \ \ \ \ \ \ \ \ \ \ \ \ \ \ \
\ \ \ \ \ \ \ \ \ \ \ \ \ \ \ \ \ 
\end{quote}

\begin{enumerate}
\item H. Amann, Linear and quasi-linear equations,1, Birkhauser, Basel 1995.

\item S. Yakubov and Ya. Yakubov, Differential-operator Equations. Ordinary
and Partial \ Differential Equations , Chapman and Hall /CRC, Boca Raton,
2000.

\item Krein S. G., Linear differential equations in Banach space, American
Mathematical Society, Providence, 1971.

\item Lunardi A., Analytic semigroups and optimal regularity in parabolic
problems, Birkhauser, 2003.

\item Dore C. and Yakubov S., Semigroup estimates and non coercive boundary
value problems, Semigroup Forum 60 (2000), 93-121.

\item Denk R., Hieber M., Pr\"{u}ss J., $R$-boundedness, Fourier multipliers
and problems of elliptic and parabolic type, Mem. Amer. Math. Soc. 166
(2003), n.788.

\item Shakhmurov V. B., Linear and nonlinear abstract equations with
parameters, Nonlinear Anal-Theor., 2010, v. 73, 2383-2397.

\item Shakhmurov V. B., V., Nonlinear abstract boundary value problems in
vector-valued function spaces and applications, Nonlinear Anal-Theor., v.
67(3) 2006, 745-762.

\item Shakhmurov V. B, Shahmurova A., Nonlinear abstract boundary value
problems atmospheric dispersion of pollutants, Nonlinear Anal-Real., v.11
(2) 2010, 932-951.

\item Shakhmurov V. B., Coercive boundary value problems for regular
degenerate differential-operator equations, J. Math. Anal. Appl., 292 ( 2),
(2004), 605-620.

\item Shakhmurov V. B., Degenerate differential operators with parameters,
Abstr. Appl. Anal., 2007, v. 2006, 1-27.

\item Shakhmurov V. B., Separable anisotropic differential operators and
applications, J. Math. Anal. Appl. 2006, 327(2), 1182-1201.

\item Shakhmurov V. B., Linear and nonlinear abstract equations with
parameters, Nonlinear Anal-Theor., 2010, v. 73, 2383-239.

\item Agarwal R., O' Regan, D., Shakhmurov V. B., Separable anisotropic
differential operators in weighted abstract spaces and applications, J.
Math. Anal. Appl. 2008, 338, 970-983.

\item Ashyralyev A., Claudio Cuevas and Piskarev S., "On well-posedness of
difference schemes for abstract elliptic problems in spaces", Numer. Func.
Anal.Opt., v. 29, No. 1-2, Jan. 2008, 43-65.

\item Weis L, Operator-valued Fourier multiplier theorems and maximal $L_{p}$
regularity, Math. Ann. 319, (2001), 735-758.\ 

\item Favini A., Shakhmurov V., Yakubov Y., Regular boundary value problems
for complete second order elliptic differential-operator equations in UMD
Banach spaces, Semigroup Forum, v. 79 (1), 2009.

\item Shahmurov R., On strong solutions of a Robin problem modeling heat
conduction in materials with corroded boundary, Nonlinear Anal. Real World
Appl., 2011,13(1), 441-451.

\item Shahmurov R., Solution of the Dirichlet and Neumann problems for a
modified Helmholtz equation in Besov spaces on an annuals, J. Differential
Equations, 2010, 249(3), 526-550.

\item Lions J. L and Peetre J., Sur une classe d'espaces d'interpolation,
Inst. Hautes Etudes Sci. Publ. Math., 19(1964), 5-68.

\item Burkholder D. L., A geometrical conditions that implies the existence
certain singular integral of Banach space-valued functions, Proc. conf.
Harmonic analysis in honor of Antonu Zigmund, Chicago, 1981,Wads Worth,
Belmont, (1983), 270-286.

\item Triebel H., \textquotedblright Interpolation theory, Function spaces,
Differential operators.\textquotedblright , North-Holland, Amsterdam, 1978.

\item Besov, O. V., P. Ilin, V. P., Nikolskii, S. M., Integral
representations of functions and embedding theorems, Nauka, Moscow, 1975.

\item Favini A., Goldstein G. R. , Goldstein J. A. and Romanelli S.,
Degenerate second order differential operators generating analytic
semigroups in $L_{p}$ and $W^{1,p}$, Math. Nachr. 238 (2002), 78 -- 102.

\item Keyantuo V., Lizama, C., Maximal regularity for a class of
integro-differential equations with infinite delay in Banach spaces, Studia
Math. 168 (2005), 25-50.
\end{enumerate}

\bigskip

\end{document}